\documentclass[letterpaper,10.5 pt,conference]{ieeeconf}
\usepackage{graphicx}
\usepackage{amssymb}
\usepackage{amsmath}
\usepackage{epstopdf}
\usepackage{algorithmic}
\usepackage{algorithm}

\usepackage{cite}

 \usepackage{tikz,pgfplots}

\IEEEoverridecommandlockouts 

\newcommand{\norm}[1]{\lVert#1\rVert}

\title{Distributed Model Predictive Consensus via the \\Alternating Direction Method of Multipliers}
\author{Tyler H. Summers$^\dag$ and John Lygeros$^\dag$%
\thanks{T. H. Summers is supported by an ETH Zurich Postdoctoral Fellowship.}%
\thanks{$^\dag$Automatic Control Laboratory, Department of Information Technology and Electrical Engineering, ETH Zurich, 8092 Zurich, Switzerland. {\tt\small \{tsummers,jlygeros\}@control.ee.ethz.ch}}
}

\begin{document}
\maketitle

\begin{abstract}
We propose a distributed optimization method for solving a distributed model predictive consensus problem. The goal is to design a distributed controller for a network of dynamical systems to optimize a coupled objective function while respecting state and input constraints. The distributed optimization method is an augmented Lagrangian method called the Alternating Direction Method of Multipliers (ADMM), which was introduced in the 1970s but has seen a recent resurgence in the context of dramatic increases in computing power and the development of widely available distributed computing platforms. The method is applied to position and velocity consensus in a network of double integrators. We find that a few tens of ADMM iterations yield closed-loop performance near what is achieved by solving the optimization problem centrally. Furthermore, the use of recent code generation techniques for solving local subproblems yields fast overall computation times. 
\end{abstract}

\section{Introduction}
One of the grand challenges facing modern society is the understanding and efficient control of large, complex technological networks, in which many dynamical subsystems interact. Such networks include power grids, the Internet, transportation networks, and unmanned vehicle formations. From a control perspective, one of the major difficulties is that there are structural constraints on information flow: each subsystem, or agent, in the network must act based on limited information about the whole network, and furthermore, even if agents could gather global information, determining optimal inputs may be computationally prohibitive.

In this paper, we consider an optimal consensus problem, in which a network of dynamically independent subsystems locally cooperate to optimize a global performance objective. There is a large and expanding literature on consensus problems; see e.g. \cite{olfati2007consensus,ren2008distributed} and references therein for an overview. However, much less attention has been given to \emph{optimally} achieving consensus in networks of dynamical systems, particularly when the subsystems are subject to state and input constraints. 

The existing literature on optimal consensus problems typically propose variations on distributed model predictive control. Model predictive control (MPC) is a well-developed and successful method for controlling systems with state and input constraints. In MPC, a model of the plant is used to predict future behavior of the system.  At each measurement instant, an optimization problem is solved to choose a sequence of inputs that optimize the predicted future behavior over a finite horizon. Only the first input in the sequence is applied to the system, and the process is repeated when new measurements are received. In distributed model predictive control schemes, each subsystem solves an optimization problem using only its own state and local state information measured or communicated from neighboring subsystems. 

Seminal work on distributed model predictive consensus includes \cite{richards2004decentralized,dunbar2006distributed,keviczky2006decentralized}. In \cite{dunbar2006distributed} and \cite{keviczky2006decentralized}, the global optimization problem is broken into smaller problems for each agent, and sufficient conditions for stability are given.  However, there is only a single communication of information between each state update rather than an iterative negotiation process, so these methods rely on neighboring agents not deviating too much from their plans and as a result are very conservative. Other related work on optimal consensus problems includes \cite{raffard2004distributed,johansson2006distributed,ferrari2007model,keviczky2008study,johansson2008decentralized,franco2008cooperative,ferrari2009model,siva2010robust}. In \cite{ferrari2009model}, consensus schemes for systems with single and double integrator dynamics are considered. Geometric properties of the optimal path followed by individual agents are used to guarantee stability. Again, there is no iterative negotiation process, so the results are very conservative. In \cite{johansson2006distributed,johansson2008decentralized,keviczky2008study}, the authors use an incremental subgradient method based on primal decomposition to negotiate an optimal consensus point. A stability analysis for these methods is given in \cite{keviczky2008study}. There has also been a significant amount of work on distributed model predictive control for dynamically coupled systems, including \cite{camponogara2002distributed,venkat2005stability,necoara2008application,stewart2010cooperative,giselsson2010distributed,doan2011iterative,doan2011dual,Conte2012,ConteSummers2012}.

Much of the work on distributed model predictive control and consensus focuses on deriving conditions that guarantee stability, but often requires extremely conservative assumptions. Moreover, while several different distributed optimization algorithms and variations thereof have been proposed to iteratively negotiate the solution to the optimization problem, the performance of these algorithms in distributed MPC problems remains poorly understood, especially in terms of the closed-loop performance of the policy. 

In this paper, we propose to use the Alternating Direction Method of Multipliers (ADMM) for solving a distributed model predictive consensus problem. ADMM is an augmented Lagrangian method originally proposed in the 1970s in \cite{glowinski1975approximation} and \cite{gabay1976dual}, and recently reviewed in \cite{boyd2011distributed}. It can be seen as a variation of dual decomposition that provides improved theoretical and practical convergence properties. The main advantages of using this method are: (1) in theory, convergence is guaranteed for any convex cost functions and constraints, and (2), in practice, the augmented Lagrangian term often speeds up convergence. We apply the method to a flocking problem in a network of double integrators and find that a few tens of ADMM iterations yield closed-loop performance near what is achieved by solving the optimization problem centrally. We also utilize recent code generation techniques for convex optimization \cite{mattingley2012cvxgen,Domahidi2012} for solving the local subproblems, thereby achieving surprisingly fast overall computation times. 

The rest of the paper is organized as follows. Section 2 formulates a distributed model predictive consensus problem. Section 3 reviews the Alternating Direction Method of Multipliers and describes how it can be used to solve the distributed model predictive consensus problem. Section 4 demonstrates the results with a numerical example. Section 5 gives concluding remarks and an outlook for future research.

\section{Distributed Model Predictive Consensus}

\subsection{An Optimal Consensus Problem}
We consider a network of $N$ independent agents, who are required to cooperatively solve an optimal control problem in a distributed way. The agents are allowed to communicate information according to a fixed undirected\footnote{Directed graphs are easily incorporated into our framework, but we restrict attention to undirected graphs to simplify the exposition.} graph $G=(V,E)$, which we refer to as the \emph{information architecture}. The vertex set $V = \{1,...,N\}$ represents the agents and the edge set $E \subseteq V\times V$ specifies pairs of agents that are permitted to communicate. If $(i,j) \in E$, we say that agents $i$ and $j$ are \emph{neighbors}, and we denote by $\mathcal{N}_i = \{j : (i,j) \in E \}$ the set of neighbors of agent $i$. 

The dynamics for the $i$th agent are given by the discrete-time linear state equation
\begin{equation}
x_i(t+1) = A_i x_i(t) + B_i u_i(t), \quad i = 1,...,N, 
\end{equation}
where $x_i(t)\in \mathbf{R}^{n_i}$ is the state, $u_i(t) \in \mathbf{R}^{m_i}$ is the control input, and $A_i \in \mathbf{R}^{n_i\times n_i}$ and $B_i \in \mathbf{R}^{n_i\times m_i}$ are given dynamics and input matrices. The global system state and input are concatenations of local states and inputs: $x(t) = [x_1(t)^T,...,x_N(t)^T]^T$ and $u(t) = [u_1(t)^T,...,u_N(t)^T]^T$.  The objective is to minimize the infinite horizon cost function
\begin{equation} \label{infhorizon}
J = \sum_{t=0}^\infty \sum_{i=1}^{N} \ell_i(x_{\mathcal{N}_i}(t),u_{\mathcal{N}_i}(t))
\end{equation}
where $\ell_i(\cdot,\cdot)$ are convex, closed, and proper stage cost functions and $x_{\mathcal{N}_i}(t)$ and $u_{\mathcal{N}_i}(t)$ are concatenations of the state and input vectors of agent $i$ and its neighbors. The stage cost functions couple the independent agents according to $G$. Each agent is subject to convex local state and input constraints $$x_{\mathcal{N}_i}(t) \in \mathcal{X}_i \quad u_{\mathcal{N}_i}(t) \in \mathcal{U}_i,$$ where $\mathcal{X}_i$ and $\mathcal{U}_i$ are convex subsets that (possibly) couple the states and inputs of neighboring agents. The global constraints are $\mathcal{X} \subseteq \mathbf{R}^{\sum_i n_i}$ and $\mathcal{U} \subseteq \mathbf{R}^{\sum_i m_i}$. In general, we seek a \emph{distributed} optimal policy $\pi : \mathcal{X} \rightarrow \mathcal{U}$ whose structure is specified by the information architecture, i.e. the policy for agent $i$ should depend only on information from its neighbors in $G$. 


\subsection{Model Predictive Control}
We employ a model predictive control scheme, in which a model is used to predict future behavior of the system, and the inputs are chosen to optimize a performance index over a finite horizon. The infinite horizon optimization problem (\ref{infhorizon}) is replaced by a finite-horizon optimization problem that is solved at each time instant in a receding horizon fashion. Only the first input of the optimizer is applied to the system; the problem is solved again whenever state information is received. The finite-horizon optimization problem is
\begin{equation} \label{mpccent}
\begin{aligned}
& \text{minimize} && J = \sum_{t=0}^{T-1} \sum_{i=1}^{N} \ell_i(x_{\mathcal{N}_i}(t),u_{\mathcal{N}_i}(t)) \\  & && \quad \quad + \sum_{i=1}^{N} \ell_{if}(x_{\mathcal{N}_i}(T),u_{\mathcal{N}_i}(T)) \\
& \text{subject to} && x_i(t+1) = A_i x_i(t) + B_i u_i(t), \\ 
&			    && x_{\mathcal{N}_i}(t) \in \mathcal{X}_i, \quad u_{\mathcal{N}_i}(t) \in \mathcal{U}_i, \\
&                             && i = 1,...,N, \quad t = 0,1,...,T-1 \\
&			    && x_{\mathcal{N}_i}(T)\in \mathcal{X}_{if} \quad i = 1,...,N,
\end{aligned}
\end{equation}
where $\ell_{if}(\cdot)$ and $\mathcal{X}_{if}$ are terminal cost functions and terminal constraint sets, respectively.
For centralized problems, there are various methods for choosing the terminal cost functions and constraint sets to ensure feasibility and stability of the closed-loop system; see e.g. \cite{mayne2000constrained}. However, for distributed problems, these quantities should be chosen to reflect the information architecture, i.e. the terminal cost and constraint sets for agent $i$ should only depend on $x_{\mathcal{N}_i}$. Recent work has considered the problem of choosing terminal costs and constraint sets with a certain distributed structure for feasibility and stability \cite{Conte2012}.

\subsection{Consistency Constraint Form}
All coupling variables in the problem can be reduced to so-called consistency constraints by introducing local copies of variables at each subsystem and requiring the local copies of coupling variables to be the same across coupled subsystems. Specifically, let $\mathbf{x}_i$ be a local variable vector for agent $i$ that includes a copy of the state and input vectors over the finite horizon and also copies of the state and input vectors over the finite horizon of all neighbors. For example, if agent $1$ is connected to agent $2$, then $\mathbf{x}_1$ includes a local copy of $x_1$, $u_1$, $x_2$, and $u_2$ over the finite horizon. Likewise, for agent $2$, the variable $\mathbf{x}_2$ contains a \emph{separate} local copy of $x_1$, $u_1$, $x_2$, and $u_2$ over the finite horizon. Let $\mathbf{x} = [\mathbf{x}_1^T,...,\mathbf{x}_N^T]^T$ be a  vector that includes all copies of all variables in the problem. The problem then has the form
\begin{equation} \label{admmified}
\begin{aligned}
& \underset{{\mathbf{x}_i \in \mathbf{X}_i}}{\text{minimize}} &&  \sum_{i=1}^{N} f_i(\mathbf{x}_i) \\
& \text{subject to} && \mathbf{x}_i = \mathbf{E}_i \mathbf{x}, \quad i=1,...,N,
\end{aligned}
\end{equation}
where $\mathbf{E}_i$ is a matrix that picks out components of $\mathbf{x}$ that match components of the local variable $\mathbf{x}_i$ (to simplify notation, there may be some redundant coupling constraints as written, but we assume that redundant coupling constraints are removed). The constraint set $\mathbf{X}_i$ is convex and includes all constraints for agent $i$ and all of its neighbors. The problem (\ref{admmified}) now has a separable cost function in variables $\mathbf{x}_i$ with only coupling equality constraints, which enforce consistency of local variable copies.

Our focus is on using a distributed optimization method to solve (\ref{admmified}) (hence, (\ref{mpccent})) by distributing the computation across a collection of communicating processors. Several different distributed optimization methods have been proposed for solving finite-horizon optimization problem arising from distributed model predictive control and consensus problems \cite{camponogara2002distributed,venkat2005stability,necoara2008application,stewart2010cooperative,giselsson2010distributed,doan2011iterative,doan2011dual,Conte2012,ConteSummers2012,richards2004decentralized,dunbar2006distributed,keviczky2006decentralized,johansson2006distributed,raffard2004distributed,ferrari2007model,keviczky2008study,johansson2008decentralized,franco2008cooperative,ferrari2009model,siva2010robust}. Many of them are based on some version of standard dual decomposition, in which a partial Lagrangian is formed by adding the coupling equality constraints to the objective function. However, these methods provably converge only under strong assumptions (e.g. that the cost function is strictly convex) and may also be slow to converge in practice. The Alternating Direction Method of Multipliers (ADMM) is an augmented Lagrangian method that guarantees convergence under very mild assumptions and also has been observed to improve convergence speed in practice. We briefly review this method in the following section. 

\section{The Alternating Direction Method of Multipliers}
ADMM was originally proposed in the 1970s in \cite{glowinski1975approximation} and \cite{gabay1976dual}, but has seen a recent resurgence in the context of increased computational power and distributed computing platforms. An excellent recent overview of the method with applications to statistics and machine learning is provided in \cite{boyd2011distributed}; the reader is referred there for details. In this section, we briefly review the highlights.

\subsection{Dual Decomposition}
Consider the optimization problem 
\begin{equation} \label{dd}
\begin{aligned}
&  \underset{{\mathbf{x}_i \in \mathbf{X}_i}}{\text{minimize}} && f(\mathbf{x}) = \sum_{i=1}^N f_i(\mathbf{x}_i) \\
& \text{subject to} && A\mathbf{x} = b
\end{aligned}
\end{equation}
where $\mathbf{x}_i \in \mathbf{R}^{n_i}$ are subvectors of the global variable $\mathbf{x} = [\mathbf{x}_1^T,...,\mathbf{x}_N^T]^T$ and $A$ is partitioned comformably $A = [A_1 \cdots A_N]$. Note that (\ref{admmified}) has the same form as problem (\ref{dd}), with $A = \text{diag}(I_{n_i}) - E$, where $E = [E_1^T,...,E_N^T]^T$, and $b=0$. The Lagrangian for problem (\ref{dd}) is \begin{equation}
\begin{aligned}
L(\mathbf{x}, \lambda) &=  \sum_{i=1}^N (f_i(\mathbf{x}_i) + \lambda^T A_i \mathbf{x}_i - \frac{1}{N} \lambda^T b) \\
		      &=  \sum_{i=1}^N L_i(\mathbf{x}_i,\lambda)
\end{aligned}
\end{equation}
where $\lambda$ are Lagrange multipliers, or dual variables, and the corresponding dual function is 
\begin{equation}
g(\lambda) = \inf_\mathbf{x} L(\mathbf{x},\lambda).
\end{equation}
Dual decomposition involves applying gradient ascent to maximize the dual function, yielding the algorithm
\begin{equation} 
\begin{aligned}
\mathbf{x}_i^{k+1} &= \underset{{\mathbf{x}_i \in \mathbf{X}_i}}{\text{argmin}} \ L_i(\mathbf{x}_i, \lambda^k) \\
\lambda^{k+1}    &= \lambda^k + \alpha^k (A \mathbf{x}^{k+1} - b)
\end{aligned}
\end{equation}
where $k$ is an iteration counter and $\alpha^k$ is a step size. Since the Lagrangian is separable in the partition of $\mathbf{x}$, the $\mathbf{x}_i$ minimizations in the first line of the algorithm can be performed independently and in parallel for all $i = 1,...,N$ for a given $\lambda^k$. The residuals $A_i \mathbf{x}_i^{k+1}-b_i$ are then gathered to update the dual variable in the second line of the algorithm. When the coupling constraint matrix $A$ has particular structure, as in problem (\ref{admmified}), the dual update need not be performed by a centralized entity, but rather can be done locally by sharing information only amongst coupled subsystems.

While dual decomposition yields a distributed algorithm for solving problem (\ref{dd}), it is only guaranteed to converge for suitably chosen step size sequences and under rather strong assumptions on the original problem. One way to improve the convergence guarantees is to augment the Lagrangian with an additional penalty term, which we discuss next.

\subsection{The Method of Multipliers}
The augmented Lagrangian includes an additional penalty term:
\begin{equation}
L_\rho(\mathbf{x}, \lambda) = f(\mathbf{x}) + \lambda^T(A\mathbf{x}-b) + \frac{\rho}{2}\norm{A\mathbf{x}-b}_2^2,
\end{equation}
where $\rho>0$ is called the penalty parameter, and can be seen as the unaugmented Lagrangian for the problem
\begin{equation} \label{mom}
\begin{aligned}
& \text{minimize} && f(\mathbf{x}) +  \frac{\rho}{2}\norm{A\mathbf{x}-b}_2^2 \\
& \text{subject to} && A\mathbf{x} = b,
\end{aligned}
\end{equation}
which has the same minimizer as problem (\ref{dd}). The penalty term has a smoothing effect that renders the dual function differentiable under mild conditions on the primal problem. The gradient ascent algorithm on the dual of (\ref{mom}) converges under very mild conditions on the primal problem. In particular, the primal objective function need not be strictly convex, finite, or even differentiable. However, while the convergence properties have improved, the quadratic penalty term prevents the $\mathbf{x}$ minimization in the dual ascent algorithm from being separable, and therefore cannot be used to distribute the computation across a set of processors, even when the objective function $f$ is separable. 

\subsection{ADMM}
ADMM can be viewed as a way to combine the robustness of the method of multipliers with the decomposability of standard dual decomposition. We now focus specifically on coupling constraints of the form given in (\ref{admmified}), which results in a particular form of ADMM that is referred to in \cite{boyd2011distributed} as ``general form global consensus''. A global variable $z$ is introduced, and the algorithm solves problems of the form
\begin{equation} \label{admm}
\begin{aligned}
&  \underset{{\mathbf{x}_i \in \mathbf{X}_i}}{\text{minimize}} && \sum_{i=1}^N f_i(\mathbf{x}_i) \\
& \text{subject to} && \mathbf{x}_i - \mathbf{\bar{E}}_i z = 0, \quad i=1,...,N.
\end{aligned}
\end{equation}
The augmented Lagrangian for problem (\ref{admm}) is
\begin{equation}
\begin{aligned}
L_\rho(\mathbf{x},z,\lambda) = &\sum_{i=1}^N \Bigl[ f_i(\mathbf{x}_i) \bigr. \\
				  &\left. + \lambda_i^T(\mathbf{x}_i - \mathbf{\bar{E}}_i z) + \frac{\rho}{2} \norm{\mathbf{x}_i - \mathbf{\bar{E}}_i z}_2^2 \right] \\
				  &= \sum_{i=1}^N L_{\rho i} (\mathbf{x}_i, z, \lambda)
\end{aligned}
\end{equation}
The ADMM algorithm consists of the iterations
\begin{equation}  \label{admmalg}
\begin{aligned}
\mathbf{x}_i^{k+1} &= \underset{{\mathbf{x}_i \in \mathbf{X}_i}}{\text{argmin}} \  L_{\rho i}(\mathbf{x}_i, z^k, \lambda^k) \\
z^{k+1} &= \underset{z}{\text{argmin}} \ L_\rho(\mathbf{x}^{k+1}, z, \lambda^k) \\
\lambda_i^{k+1}    &= \lambda_i^k + \rho (\mathbf{x}_i^{k+1} - \mathbf{\bar{E}}_i z^{k+1} )
\end{aligned}
\end{equation}
The method of multipliers for problem (\ref{admm}) would jointly minimize over $\mathbf{x}$ and $z$, whereas here the algorithm alternates sequentially between $\mathbf{x}$ and $z$  minimizations, which allows the $\mathbf{x}_i$ minimizations to be done in parallel. In this particular form, the $z$ minimization also distributes across its components and can be simplified to
\begin{equation}
z_i^{k+1} = \frac{1}{\vert \mathcal{N}_i \vert + 1 } \sum_{j \in \mathcal{N}_i} \left[(\mathbf{x}_j^k)_i + \frac{1}{\rho}(\lambda_j^k)_i \right],
\end{equation}
where $(\mathbf{x}_j^k)_i$ and $(\lambda_j^k)_i$ denote the components of $\mathbf{x}_j$ and $\lambda_j$ that map to the $i$th component of $z$. 
It can be shown that after one iteration, the local average of the dual variables for each component of $z$, viz. $\sum_{j \in \mathcal{N}_i} \frac{1}{\rho}(\lambda_j^k)_i$, is zero. Hence, the $z$ update further reduces to
\begin{equation}
z_i^{k+1} = \frac{1}{\vert \mathcal{N}_i \vert + 1 } \sum_{j \in \mathcal{N}_i \cup i} (\mathbf{x}_j^k)_i,
\end{equation}
i.e. the $z$ update is simply a local averaging of all ``opinions'' of each agent on the particular component of $z$ (see \cite{boyd2011distributed}). 

\begin{algorithm}[h!]
\caption{Alternating Direction Method of Multipliers}
\label{alg:admm}
\begin{algorithmic}[1]
\STATE $\forall i$ in parallel:
\STATE initialize $\lambda = 0$, $z = 0$
\REPEAT
\STATE $\mathbf{x}_i^{k+1} = \underset{\mathbf{x}_i \in \mathbf{X}_i}{\arg\min} \  L_{\rho}^i(\mathbf{x}_i,z^k,\lambda^k)$
\STATE communicate $\mathbf{x}_i^+$ to all $j \in \mathcal{N}_i$
\STATE $z_i^{k+1} = \frac{1}{\vert \mathcal{N}_i \vert + 1 } \sum_{j \in \mathcal{N}_i \cup i} \ (\mathbf{x}_j^k)_i$
\label{algstep:admm_average}
\STATE $\lambda_i^{k+1}   = \lambda_i^k + \rho (\mathbf{x}_i^{k+1} - \mathbf{\bar{E}}_i z^{k+1} )$
\label{algstep:admm_lambda_update}
\UNTIL{convergence or maximum iterations reached}
\end{algorithmic}
\end{algorithm}

ADMM is a distributed message passing algorithm. In the context of the distributed model predictive control problem, the algorithm works as follows. Each agent solves a local subproblem to minimize its local cost function, which computes a local plan for both its own state and input and also for states and inputs of neighboring agents over the horizon. The solution of this problem (step 1 in (\ref{admmalg})) can be interpreted as the agent's own ``opinion'' on what it and its neighbors should do to minimize its local cost. Neighboring agents then communicate their plans to one another. The global variable is updated (step 2 in (\ref{admmalg})) by averaging the local plans. The dual variables are then updated (step 3 in (\ref{admmalg})) based on the difference between the local plans and the averaged global plan. The whole iterative negotiation process takes place within each sampling period, every time a state measurement is obtained.

At each iteration of the algorithm, each subsystem solves a small (relative to the global problem), structured convex optimization problem. Very recently, code generation techniques for convex optimization have been developed \cite{mattingley2012cvxgen,Domahidi2012}. By exploiting problem structure, these techniques can enable extremely fast computation times for small convex optimization problems (on the order of microseconds to milliseconds, depending on problem size and processor capabilities). ADMM breaks a large optimization problem into many smaller optimization problems that can be spread across multiple processors. Combining code generation techniques with a distributed optimization method like ADMM potentially allows solving large optimization problems surprisingly fast.

The distributed model predictive control literature often lacks evaluations of proposed distributed optimization algorithms in terms of computation times and performance comparison to centralized schemes. In the next section, we explore the practical performance of ADMM on a distributed model predictive consensus problem.

\section{Numerical Example}
In this section, we apply the method to a consensus problem for a  collection of double integrators moving in three-dimensional space with input constraints. The state of each agent $x_i \in \mathbf{R}^6$ consists of spatial positions (the first, third, and fifth components) and velocities (the second, fourth, and sixth components). The input $u_i \in \mathbf{R}^3$ consists of force components acting on the agent. After an Euler discretization with step size $T_s$, the dynamics and input matrices for each agent are
\begin{equation}
A_i = I_3 \otimes \left[\begin{array}{cc}1 & T_s \\0 & 1\end{array}\right], \quad B_i = I_3 \otimes \left[\begin{array}{c}0 \\  T_s/ m_i\end{array}\right].
\end{equation}
where $\otimes$ is the Kronecker product and $m_i$ is the mass of agent $i$. We consider input constraints only: $\mathcal{X}_i = \mathbf{R}^6$ and $\mathcal{U}_i = \{ u \in \mathbf{R}^3 \ \vert \ \norm{u}_\infty \leq u_{max} \}$.

We consider a flocking problem, in which the objective is for all agents to converge to a common position and velocity. One possible stage cost function that encodes this is
\begin{equation}
\ell(x(t),u(t)) = x(t)^T (I_3 \otimes (L \otimes I_2)) x(t) + u(t)^T R u(t),
\end{equation}
where $L$ is the \emph{weighted graph Laplacian matrix}\footnote{The weighted graph Laplacian matrix $L \in \mathbf{R}^{N\times N}$ for a graph $G=(V,E)$ is defined as follows. We first define the weighted adjacency matrix $A \in \mathbf{R}^{N \times N}$, which encodes the communicated state information flow. The diagonal entries of $A$ are zero, and an off-diagonal entry $a_{ij}$ is positive if agents $i$ and $j$ are neighbors. Then the Laplacian matrix is defined by $l_{ii} = \sum_i a_{ij} $ and $l_{ij} = -a_{ij}$. } and $R$ is a positive-definite matrix. In the simulations in the section, we use $R = I_3 \otimes I_N$. We also include process noise that acts independently on each component of the acceleration of each agent and has zero-mean and covariance $0.1$. We consider $N=5$ agents connected in a simple path graph. The prediction horizon is $T=10$. 

To evaluate the method, we generated random initial conditions and random noise sequences and ran a closed-loop simulation for 250 steps. At each step, the agents obtain a state measurement and communicate their state to their neighbors. Then the agents commence the iterative ADMM algorithm described in the previous section, communicating the intermediate results of local optimizations to neighbors and locally updating global variables and Lagrange multipliers. We fixed the maximum number of ADMM iterations to 30 and compared the closed-loop performance with the performance achieved when solving the finite-horizon model predictive control problem centrally.

Figures 1-3 show example state and input time trajectories for each spatial component for a randomly generated initial condition and noise sequence. The solid blue curve shows the trajectory using ADMM in the closed-loop, and the dashed red curve shows the trajectory with the finite horizon MPC problem solved with a centralized solver. Figure 4 shows a spatial plot of the 5 agents moving through space. We can see that both methods achieve consensus, but converge to slightly different points on the consensus manifold. Figure 5 shows a comparison of closed-loop performance (expressed as percentage of total cost) achieved by a centralized solver with the performance of ADMM with various maximum number of iterations. For each maximum iteration number, we randomly generated 120 initial conditions and noise sequences, ran a 250-step closed-loop simulation,  and averaged the results. If only 1 iteration is performed, the performance is about 80\% worse; however, for this example, using just 2 iterations achieves performance within 1.5\% of the centralized solver and using 10 or more achieves performance within 0.5\% (on average).

We used C code generated from CVXGEN \cite{mattingley2012cvxgen} to solve the local subproblems. Using a 3.1 GHz Intel Core i5 processor, the local subproblems are solved in under 2 milliseconds; thus, if the local subproblems were solved in parallel, the total computation time to finish 30 negotiations to compute the optimal inputs at each sampling time would be under 60 milliseconds. There is evidence that ADMM scales very well with problem size; a computational study in the context of distributed model predictive control is performed in \cite{ConteSummers2012}, and in \cite{kraning2012message}, a dynamic energy network scheduling problem with 10 million variables was reported to be capable of being solved in a parallel implementation in under 1 second. Thus, there is potential to solve very large problems surprisingly fast by combining  distributed optimization methods like ADMM with fast custom solvers for local subproblems.

\begin{figure}
\begin{center}
\resizebox{0.94\linewidth}{!}{\includegraphics{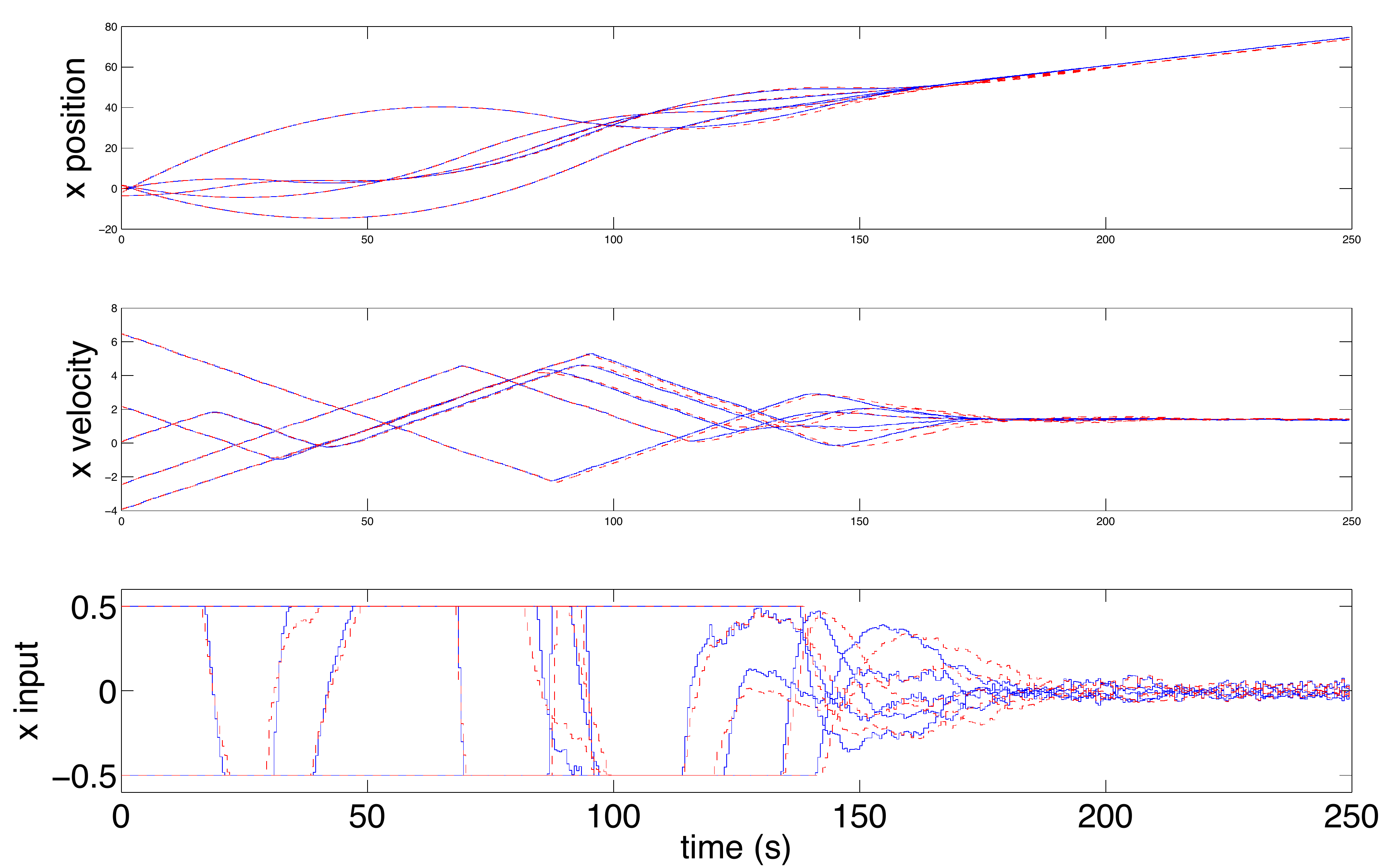}}
\end{center}
\vspace{-0.5cm}
\caption{$x$ component time trajectories.}
\label{fig1}
\vspace{-0.4cm}
\end{figure}

\begin{figure}
\begin{center}
\resizebox{1.06\linewidth}{!}{\includegraphics{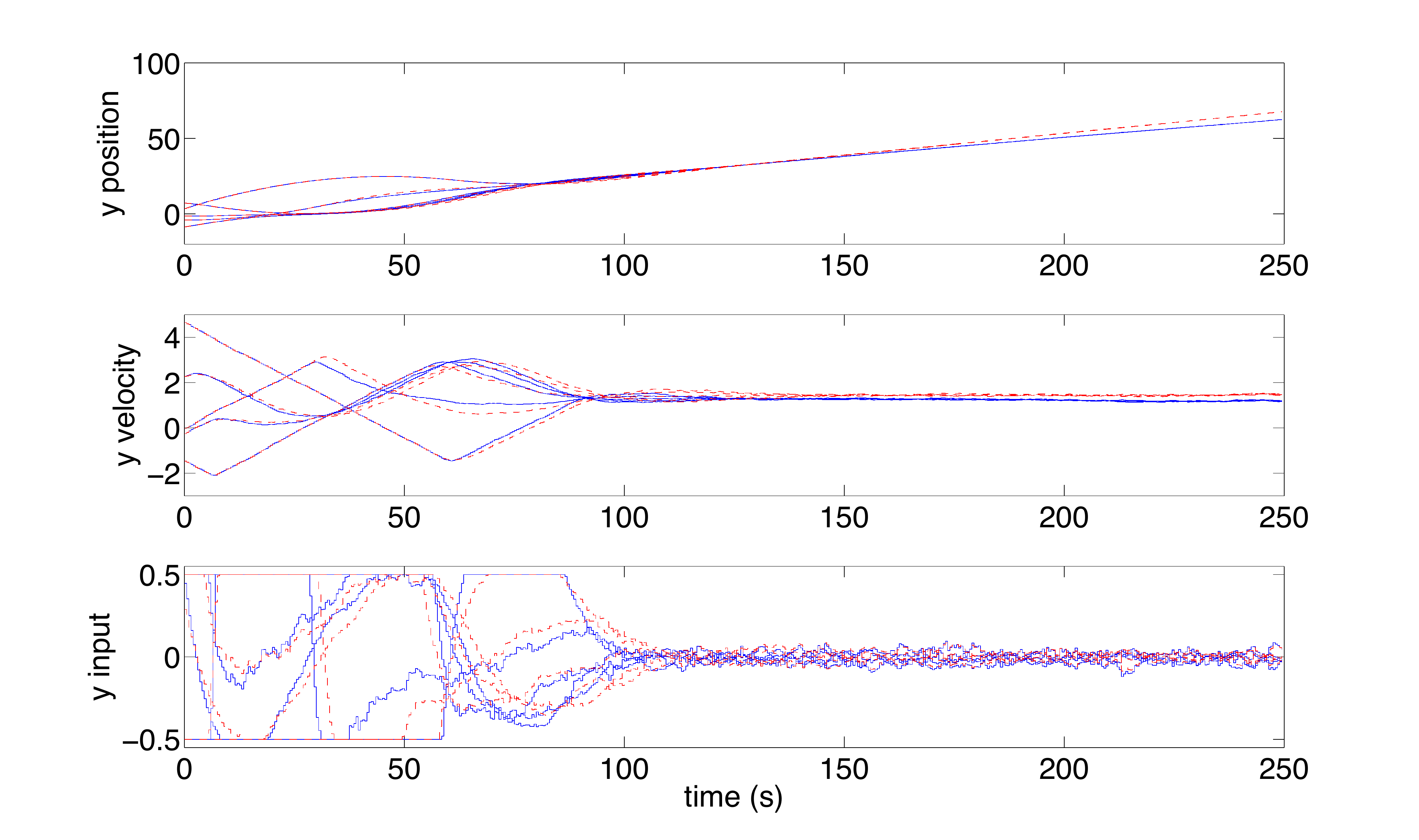}}
\end{center}
\vspace{-0.5cm}
\caption{$y$ component time trajectories.}
\label{fig2}
\vspace{-0.7cm}
\end{figure}

\begin{figure}
\vspace{-0.7cm}
\begin{center}
\resizebox{1.06\linewidth}{!}{\includegraphics{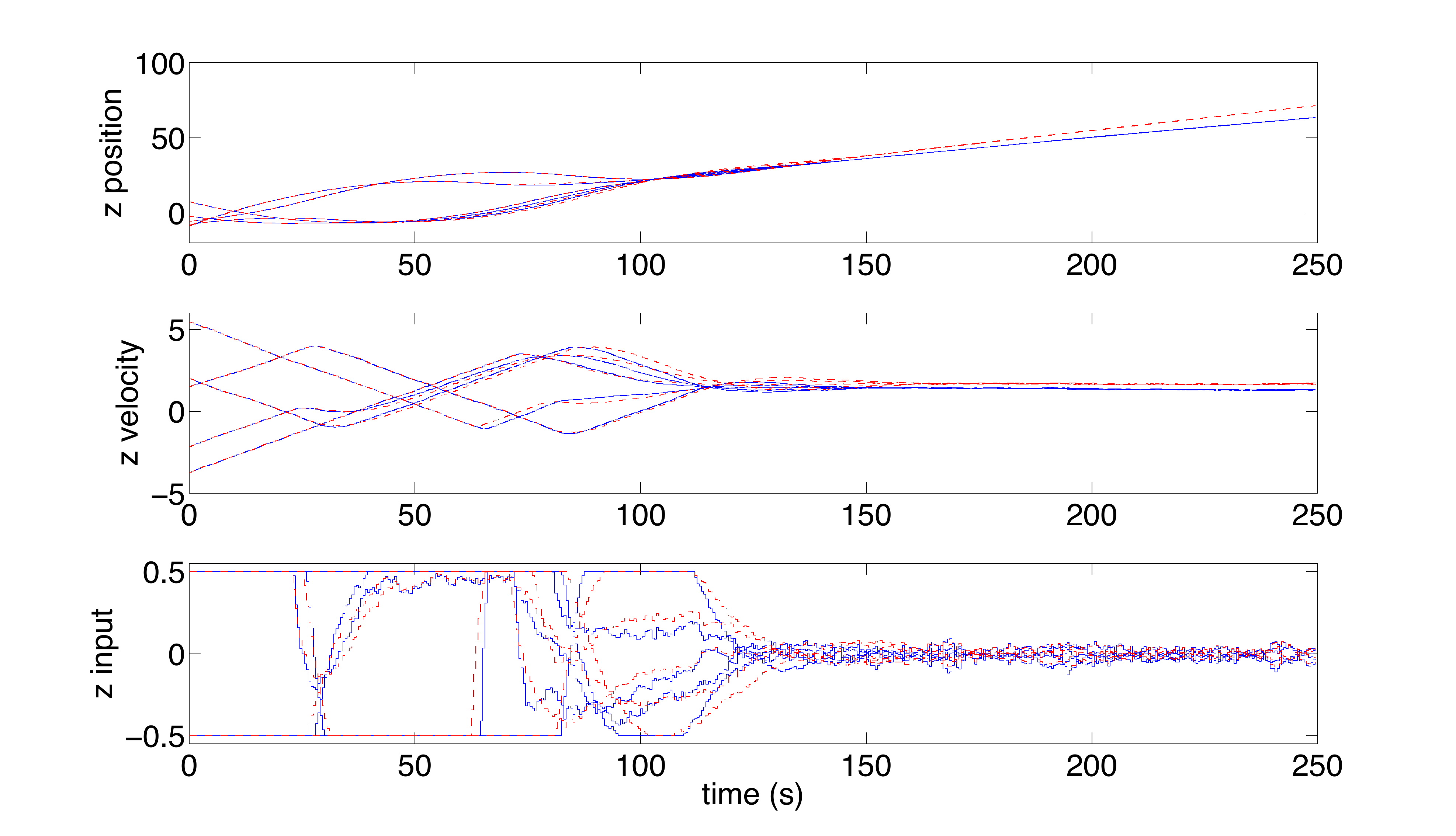}}
\end{center}
\vspace{-0.5cm}
\caption{$z$ component time trajectories.}
\label{fig3}

\end{figure}

\begin{figure}
\begin{center}
\resizebox{1.04\linewidth}{!}{\includegraphics{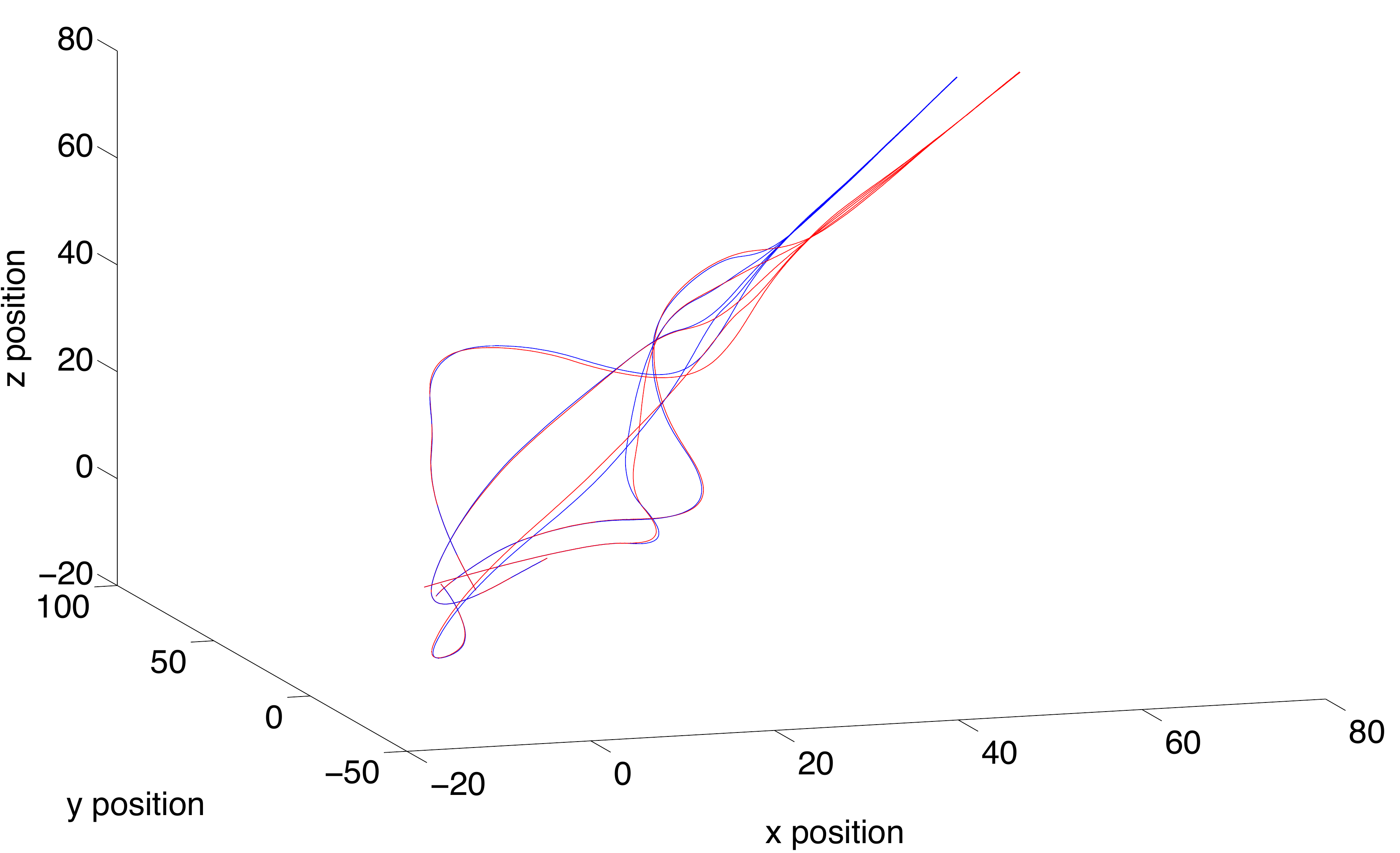}}
\end{center}
\vspace{-0.5cm}
\caption{Spatial plot of agent trajectories.}
\label{fig4}
\vspace{-0.4cm}
\end{figure}

\begin{figure}
\begin{center}
\resizebox{0.86\linewidth}{!}{\includegraphics{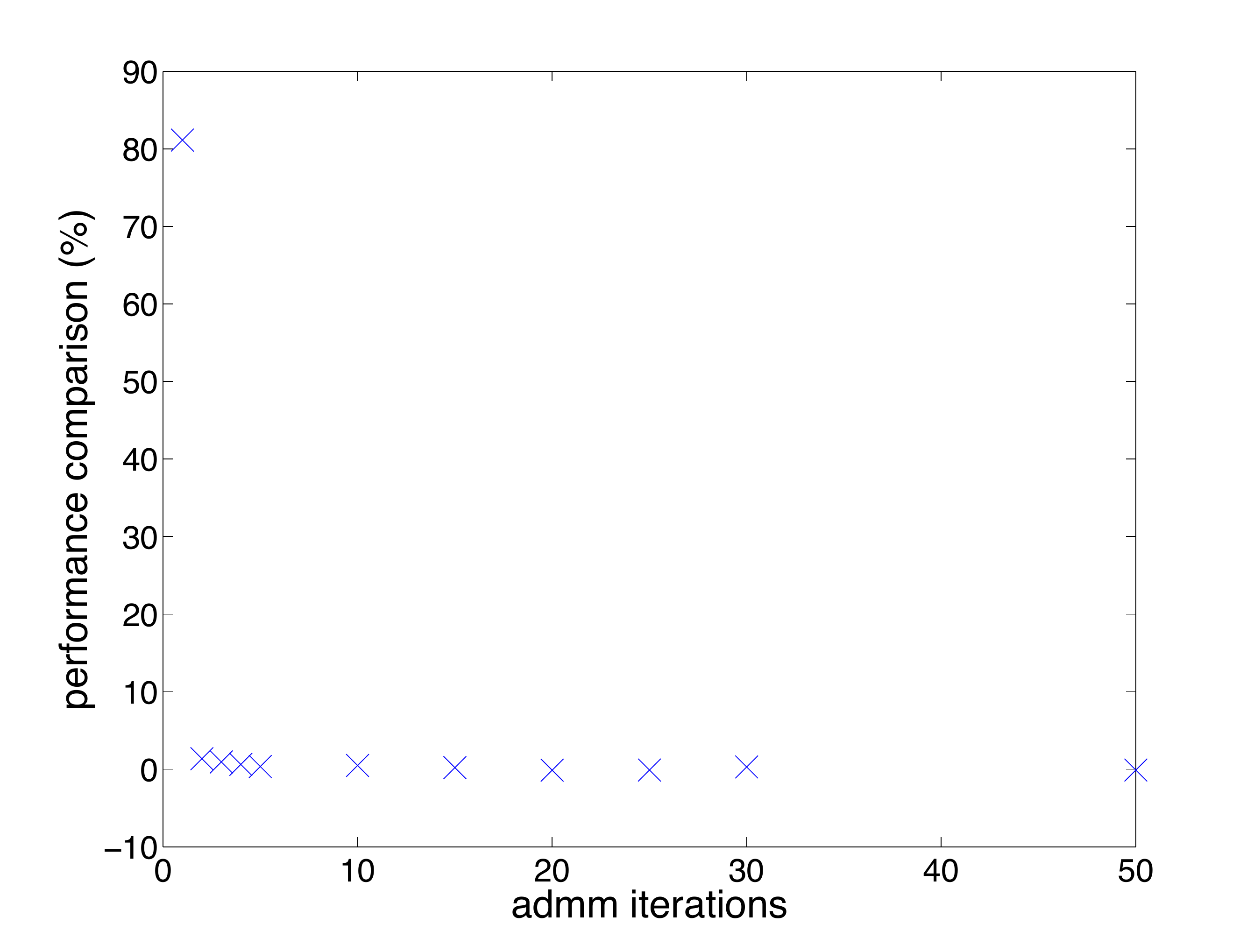}}
\end{center}
\vspace{-0.5cm}
\caption{Closed-loop performance comparison.}
\label{fig5}
\vspace{-0.6cm}
\end{figure}


\section{Conclusions}
In this paper, we proposed the Alternating Direction Method of Multipliers (ADMM), which has good theoretical and practical convergence properties, for solving a distributed model predictive consensus problem. A numerical example illustrates that near-centralized performance can be achieved with only a few tens of iterations of the distributed method.

There are many potential directions for future research. On the theoretical side, while ADMM has good theoretical and practical convergence properties, the consistency constraints of the iterates are only feasible in the limit. So guaranteeing closed-loop stability for early stopping is an open problem and will require arguments different from the standard ones used in model predictive control. Communication network non-idealities, e.g. delays and packet drops, are important and may dominate computation times in a distributed implementation, and should be considered in future research. Finally, there are many potential applications for distributed model predictive control and consensus, e.g. to vehicle formations and power networks.
\vspace{-0.15cm}

\bibliographystyle{IEEETran}  
\bibliography{refs}  

\end{document}